\documentclass[a4paper, 11pt]{article}
\usepackage[english]{babel}
\usepackage{amsmath, amssymb, amsfonts, amsthm, bm}
\usepackage{mathrsfs}
\usepackage[mathcal]{euscript}
\usepackage{enumerate}
\usepackage{mathbbol}
\usepackage{float}
\usepackage{color,graphicx,epstopdf}
\usepackage{dsfont}
\usepackage{hyperref}
\usepackage{caption}
\usepackage{subcaption}
\usepackage{cite} 
\usepackage{algorithm} 
\usepackage{algorithmic} 
\usepackage{tikz-cd}
\usepackage{mathtools}

\DeclareMathOperator*{\argmin}{arg\,min}

\DeclareFontFamily{OT1}{pzc}{}
\DeclareFontShape{OT1}{pzc}{m}{it}{<-> s * [1.200] pzcmi7t}{}
\DeclareMathAlphabet{\mathpzc}{OT1}{pzc}{m}{it}

\newcommand*\mcapinn[2]{\vcenter{\hbox{$\mathsurround=0pt
  \ifx\displaystyle#1\textstyle\else#1\fi\bigcap$}}}

\newcommand*\mcupinn[2]{\vcenter{\hbox{$\mathsurround=0pt
  \ifx\displaystyle#1\textstyle\else#1\fi\bigcup$}}}

\newtheorem{theorem}{Theorem}
\newtheorem{definition}{Definition}

\newtheorem{proposition}{Proposition}

\title{\bf Social Shaping of Competitive Equilibriums\\ for Resilient   Multi-Agent Systems}
\date{}

\author{Yijun Chen, Razibul Islam, Elizabeth Ratnam, Ian R. Petersen, and Guodong Shi\thanks{Y. Chen and G. Shi are with the  Australian Center for Field Robotics, The University of Sydney, NSW, Australia.
(E-mail: yijun.chen@sydney.edu.au; guodong.shi@sydney.edu.au)} \thanks{R. Islam, E. Ratnam and I. R. Petersen are with the Research School of Engineering, The Australian National University, Canberra, Australia. (E-mail:  sk.islam@anu.edu.au; elizabeth.ratnam@anu.edu.au; ian.petersen@anu.edu.au) }
} 
\begin{document}
\maketitle

\begin{abstract}
In this paper, we study multi-agent systems with decentralized resource allocations. Agents have local demand and resource supply, and are interconnected through a network designed to support sharing of the local resource; and the network has no external resource supply. It is known from classical welfare economics theory that by pricing the flow of resource, balance between the demand and supply is possible. Agents decide on the consumed resource, and perhaps further the traded resource as well, to maximize their payoffs considering both the utility of the consumption, and the income from the trading. When the network supply and demand are balanced, a competitive equilibrium is achieved if all agents maximize their individual payoffs, and a social welfare equilibrium is achieved if the total agent utilities are maximized. First, we consider multi-agent systems with static local allocations, and prove from duality theory that under general convexity assumptions, the competitive equilibrium and the social welfare equilibrium exist and agree. Compared to similar results in the literature based on KKT arguments, duality theory provides a direct way for connecting the two notions and for a more general (e.g. nonsmooth) class of utility functions. Next, we show that the agent utility functions can be prescribed in a family of socially admissible functions, under which the resource price at the competitive equilibrium is kept below a threshold. Finally, we extend the study to dynamical multi-agent systems where agents are associated with dynamical states from linear processes, and we prove that the dynamic the competitive equilibrium and   social welfare equilibrium continue to exist and coincide with each other.

\end{abstract} 

\section{Introduction}

Next generation technologies are leveraging the internet of things (IoT) to support critical infrastructure systems including energy distribution and automotive transportation, and are being organized as interconnected multi-agent systems \cite{Mesbahi2010}. Such systems involve  data collection, resource allocation, and control  coordination between geographically distributed subsystems. Each subsystem, termed an `agent', is an intelligent functioning unit with its own decision, objective and preference, and   remarkably,  network-level goals such as consensus, formation, and optimality can be achieved by agents interacting with others over a {\em network} e.g., \cite{tsi1984,Jadbabaie03,Murray04,Martinez07,Nedic10}. The underlying network for multi-agent systems can be physical such as transmission lines in a power grid, non-physical such as wireless communication channels, or a combination of the two.  The key promise of organizing subsystems into networked  multi-agent systems is a radical improvement in scalability, efficiency,   and sustainability through shared inputs and outputs, and coordinated decisions and controls. 

One important problem for multi-agent system operation is efficient  resource allocation, where demand and supply must be balanced for efficient and secure operations at the system level.   
In a typical resource allocation problem, agents have local demand and internal and external resource  suppliers,    interconnected through a network  that allows for transmission of the  resource.  In light of classical welfare economics theory \cite{Debreu1952,Arrow1954}, careful  pricing of the transmission flow potentially balances the demand against the supply across the entire system.  Agents decide on the resource consumed, and perhaps further the resource traded, to maximize  their payoffs considering both the utility from production/consumption, and income from the trading. When network supply and demand is balanced, a competitive    equilibrium is achieved if all agents maximize  their  individual payoffs; a social welfare equilibrium is achieved if the total agent utilities are maximized \cite{Acemoglu}. 

The concept of resource allocation via a competitive equilibrium has been widely applied in the smart grid literature, where consumers and suppliers of electrical energy make trading decisions over an energy market to achieve a competitive equilibrium  \cite{Chen10,Li2011,Papadaskalopoulos13,Hansen15,Li2015,
Jadhav2018,Muthirayan2019}. The dynamics of power networks can even be coupled with the market pricing dynamics \cite{Alvarado2001,Singh2018,Stegink2017,Zhang2015,Jovic2010}, where energy price becomes an effective controller for quasi steady-state grid frequency stability. We refer to \cite{CSM2020} for a comprehensive survey on transactive energy systems. Moreover,  in climate-economy  frameworks, a   market approach    was also proposed as a principled way for mitigating global carbon emissions while balancing regional interests  \cite{Nordhaus96,Nordhaus2011}. The price for  carbon  emissions is calculated under the competitive equilibrium of the carbon market, becoming a benchmark for the social cost of carbon \cite{Nordhaus2017}. We refer to  \cite{Chris2019} for an excellent introduction  into the dynamic integration of climate and economy models from a feedback system perspective.

Despite the aforementioned successes in coordinating   multi-agent systems via market pricing in critical engineering and societal problems, the resilience of the pricing mechanism is potentially a serious challenge even for theoretically optimal equilibrium conditions. In fact, it is well known that in welfare economics that an efficient market equilibrium may imply  Pareto optimality, while individual  equity and fairness of the agents may  be discarded entirely \cite{Hallam}. In the context of power grids, market-driven prices have exceeded affordable thresholds for an individual household. Specifically, in February 2021, Texas encountered extreme cold weather conditions resulting in a power outage disaster throughout the state, with customers on rolling blackouts.  In Texas, competitive pricing did not curtail non-critical loads to balance the grid under the power shortage conditions \cite{CNN}. Residential customers in Texas paying wholesale electricity prices during the power shortage event reported electricity bill shock with their bills exceeding previous invoices by a factor of over one hundred. Moreover, in the context of climate change mitigation, it was only possible for the Paris agreement to be reached in 2016 after a decade-long negotiation, as the  estimated social cost of carbon was perceived as unfair among different nations  \cite{Kuyper}. It follows that a competitive equilibrium for multi-agent systems are in need of social shaping at the network level: only equilibriums within a prescribed range of fairness can be accepted; only individual agent utility functions within a prescribed range of social responsibilities can be admissible. More importantly, the two directions of social shaping should be consistent: socially admissible agent utility functions should always lead to socially acceptable equilibriums.

In this paper, multi-agent systems with decentralized resource allocations are entirely self-sustained, i.e. there is no external resource supply.  First, we consider multi-agent systems with static local allocations, and we prove   that under general convexity assumptions, the competitive equilibrium and the social welfare equilibrium  exist and agree. Our proof is based on duality theory, as opposed to the commonly attributed KKT arguments (e.g. \cite{Li2015}) that support similar results. Here duality theory provides  a  more direct way of  connecting the competitive equilibrium to the social welfare equilibrium, and supports a more general (e.g. nonsmooth) class of utility functions. Next, we investigate the case when the pricing under a competitive equilibrium is associated with an upper bound for social acceptance. By means of constructive analysis, we show that the agent utility function is prescribed by a family of socially admissible quadratic  functions, under which the pricing at the competitive equilibrium is always below a threshold. Finally, we extend the study to dynamical multi-agent systems where agents are associated with dynamical states from linear processes, and prove that the dynamic  competitive equilibrium and  the social welfare equilibrium continue to exist and coincide in an optimal control context. 

The remainder of the paper is organized as follows. In Section~\ref{sec:static}, we introduce the multi-agent system with static decisions. In Section~\ref{sec:shaping}, we formulate social shaping of competitive equilibriums. In Section~\ref{sec:dynamic}, we formulate dynamic pricing for resource allocation of multi-agent systems with an underlying dynamical process. Numerical examples are presented in Section~\ref{sec:numerical} and concluding remarks are presented in Section~\ref{sec::conc}.

\section{Static Multi-Agent Systems}\label{sec:static}

\subsection{Competitive Equilibrium for Static  Multi-Agent Systems}  
We consider a multi-agent system (MAS) with $n$ agents. The agents are  indexed in $\mathrm{V}=\{1,\dots,n\}$. We consider a basic MAS setup  with static agent decisions on load allocations. 

\medskip

\noindent {\bf MAS with Static Agent Load Decisions} (MAS-SALD).  Each agent $i$ holds a local resource of $a_i$ units, and make a (static) decision to allocate $x_i\in \mathbb{R}^{\geq 0}$ units of load for itself. The utility function related to agent $i$ allocating $x_i$ amount of load is $f_i(x_i):\mathbb{R}^{\geq 0} \mapsto \mathbb{R}$.    Consequently, agent $i$ would incur an $a_i-x_i$ amount of surplus ($a_i>x_i$), or a shortcoming ($a_i<x_i$). We assume that there is a connected   network among the $n$ agents so that they can balance the surplus and  shortcomings  through a pricing mechanism. To be precise, each unit of resource across the network is priced at $\lambda\in \mathbb{R} $. Therefore, agent $i$ will yield $(a_i-x_i)\lambda$ in income or spending.

\medskip

Denoting $\mathbf{x}=(x_1\dots x_n)^\top$ as the network resource allocation profile, we introduce the following definitions. 
\begin{definition}
A pair of price-allocation decisions ($\lambda^\ast$, $\mathbf{x}^\ast$) is a competitive equilibrium for the MAS-SALD if the following conditions hold:

(i) each agent $i$ maximizes her combined payoff at ${x}^\ast_i$, i.e.,  ${x}^\ast_i$ is an optimizer for the following constrained optimization problem: 
\begin{equation}
\begin{aligned}
\max_{{x}_i} \quad &  f_i(x_i) + \lambda^\ast (a_i -x_i) \\
{\rm s.t.} \quad & x_i\in \mathbb{R}^{\geq0}.
\end{aligned}\label{load_competitive_eq}
\end{equation}

(ii) the total demand and supply are balanced across the network: 
\begin{equation}\label{load_demand_supply_constraints}
	\sum_{i=1}^n x_i^\ast =\sum_{i=1}^n a_i.
\end{equation}
\end{definition}

\begin{definition}
A resource allocation profile $\mathbf{x}^\star$ is a social welfare  equilibrium for the MAS-SALD if it is a solution to the following optimization problem:  
\begin{equation}
\begin{aligned}
\max_{\mathbf{x}} \quad & \sum_{i=1}^n f_i(x_i) \\
{\rm s.t.} \quad & \sum_{i=1}^n x_i =\sum_{i=1}^n a_i\\
\quad & x_i\in \mathbb{R}^{\geq 0}; \; i\in\mathrm{V}. 
\end{aligned}\label{eq1}
\end{equation}
\end{definition}

We present the following result which establishes the equivalence between a competitive equilibrium and a social welfare equilibrium. The result is based only on a concavity assumption for the utility functions $f_i$. 
\begin{theorem}\label{thm1}
Consider the MAS-SALD. Suppose each $f_i(\cdot)$ is   concave over the domain $\mathbb{R}^{\geq 0}$. Then the social welfare  equilibrium(s)  and the competitive equilibrium(s)  coincide. To be precise, the following statements hold.

(i) If  ($\lambda^\ast$, $\mathbf{x}^\ast$) is a competitive equilibrium, then $\mathbf{x}^\ast$ is a social welfare  equilibrium. 

(ii) If $\mathbf{x}^\star$ is a social welfare  equilibrium, then there exists $\lambda^\ast \in \mathbb{R}$ such that ($\lambda^\ast$, $\mathbf{x}^\star$) is a competitive equilibrium.
\end{theorem}
{\it Proof.} (i) Let ($\lambda^\ast$, $\mathbf{x}^\ast$) be a competitive equilibrium. The proof proceeds by contradiction. Suppose that $\mathbf{x}^\ast$ is not a social welfare  equilibrium. Then there must exist $\bar{\mathbf{x}}^\ast$ such that $\sum_{i=1}^n \bar{x}_i^\ast =\sum_{i=1}^n {x}_i^\ast =\sum_{i=1}^n a_i$, and $\sum_{i=1}^n f_i({x}_i^\ast)<\sum_{i=1}^n f_i(\bar{x}_i^\ast)$. Consequently, there holds
\begin{align}
\sum_{i=1}^n \Big( f_i({x}_i^\ast)+ \lambda^\ast (a_i-{x}_i^\ast ) \Big) <\sum_{i=1}^n \Big( f_i(\bar{x}_i^\ast)+\lambda^\ast (a_i-\bar{x}_i^\ast ) \Big). 
\end{align}
This implies that there is at least one $m\in\mathrm{V}$ such that 
$$
f_m({x}_m^\ast)+ \lambda^\ast (a_m-{x}_m^\ast )< f_m(\bar{x}_m^\ast)+\lambda^\ast (a_m-\bar{x}_m^\ast),
$$
which contradicts the fact that ($\lambda^\ast$, $\mathbf{x}^\ast$) is a competitive equilibrium. 

(ii) We pretend a proof using duality. To be consistent with the literature on duality theory for continuous optimization, we  denote $g_i=-f_i$, and rewrite (\ref{eq1}) as 
\begin{equation}
\begin{aligned}
\min_{\mathbf{x}} \quad & \sum_{i=1}^n g_i(x_i) \\
{\rm s.t.} \quad & \sum_{i=1}^n x_i =\sum_{i=1}^n a_i; x_i\in \mathbb{R}^{\geq 0}, i=1,\dots,n .
\end{aligned}\label{eq3}
\end{equation}
Let $\mathbf{x}^\star$ be a social welfare  equilibrium\footnote{Note that $\mathbf{x}^\star$ must be finite as the feasible set of $\mathbf{x}$ is compact. }. Then from its definition there holds $\sum_{i=1}^n  {x}_i^\star =\sum_{i=1}^n a_i$. Since (\ref{eq3}) is a convex optimization problem with a linear equality constraint, strong duality holds \cite{boyd} and we denote the optimal primal and dual costs of (\ref{eq3}) as $p_\ast$ and $d^\ast$, respectively.

The Lagrangian function of (\ref{eq3}) is 
$$
L(\mathbf{x}, \lambda)=\sum_{i=1}^n g_i(x_i) + \lambda\Big( \sum_{i=1}^n x_i-\sum_{i=1}^n a_i\Big): (\mathbb{R}^{\geq 0})^{n}\times \mathbb{R} \mapsto \mathbb{R}. 
$$
Then we introduce
$$
L^\ast(\lambda)=\min_{\mathbf{x}\in (\mathbb{R}^{\geq 0})^{n}} L(\mathbf{x}, \lambda).
$$
If $\lambda^\ast$ is dual optimal (i.e., $\lambda^\ast \in \arg\max_{\lambda\in \mathbb{R}} L^\ast(\lambda)$), there holds from strong duality \cite{boyd} that
\begin{align}
d^\ast=L^\ast(\lambda^\ast)&=\min_{\mathbf{x}\in (\mathbb{R}^{\geq 0})^{n}} L(\mathbf{x}, \lambda^\ast)\\
&\leq L(\mathbf{x}^\star, \lambda^\ast)\\
&= \sum_{i=1}^n g_i({x}_i^\star)\\
&= p^\ast. 
\end{align}
This implies the inequality from the above equation actually holds at equality: 
\begin{align}\label{ee}
\mathbf{x}^\star\in\arg\min_{\mathbf{x}\in (\mathbb{R}^{\geq 0})^{n}} L(\mathbf{x}, \lambda^\ast).
\end{align}
 Note that $L(\mathbf{x}, \lambda^\ast)=\sum_{i=1}^n \big(g_i(x_i) + \lambda^\ast(x_i- a_i)\big)$ implies
 \begin{align}\label{ss}
 {x}^\star_i \in \arg \max_{x_i\in \mathbb{R}^{\geq 0}}\big( f_i(x_i) + \lambda^\ast (a_i -x_i)\big). 
  \end{align}
Thus, we have proved that $(\lambda^\ast, \mathbf{x}_i^\star)$ is a competitive equilibrium. \hfill$\square$

Clearly, in this basic multi-agent system setup, the price $\lambda^\ast$ associated with a competitive equilibrium could take negative values. From an economic point of view,  the resource at every agent must either be consumed or traded, and in cases of an oversupply of resource a negative price for load balancing would occur. From an optimization point of view, the price $\lambda^\ast$ is the Lagrangian multiplier associated with an equality constraint for a constrained optimization problem, which can take positive or negative signs. The following result indicates that as long as one agent is associated with a non-decreasing utility function, oversupply will not happen.

\begin{proposition}\label{prop1}Consider the MAS-SALD.
Suppose each $f_i(\cdot)$ is   concave over the domain $\mathbb{R}^{\geq 0}$. Let ($\lambda^\ast$, $\mathbf{x}^\ast$) be a competitive equilibrium.  Then $\lambda^\ast\geq 0$ if there exists at least one agent $m\in\mathrm{V}$ such that $f_m(\cdot)$ is non-decreasing. 
\end{proposition}
{\it Proof.} Let $f_m(\cdot)$ be non-decreasing. Assume $\lambda^\ast<0$. Then $f_m(x_m) + \lambda^\ast (a_m -x_m\big)$ is a strictly increasing function with respect to $x_m$. Therefore, there can not be a finite ${x}^\ast_m$ such that  ${x}^\ast_m\in \arg \max_{x_m\in \mathbb{R}^{\geq 0}}\big( f_m(x_m) + \lambda^\ast (a_m-x_m)\big)$, contradicting the definition of the competitive equilibrium.  \hfill$\square$

\subsection{MAS  with Trading  Decisions}
In our standing multi-agent system model, agents only decide on their allocated load $x_i$, and the surplus/shortcoming $a_i-x_i$ have to go to the network. Now we relax this restriction, and introduce the following extended MAS.

\medskip

\noindent{\bf MAS with Static Agent Load and Trading  Decisions} (MAS-SALTD) On top of the MAS-SALD,  each agent $i$   further makes a decision on the traded amount of resource, denoted $e_i$. This $e_i$   is physically constrained by $x_i$ and $a_i$ in the following way\begin{itemize}
\item[(i)] if $x_i<a_i$, then agent $i$ can sell, in which case $e_i \geq 0$ and $e_i\leq a_i-x_i$;
\item[(ii)] if $x_i\geq a_i$, then agent $i$ can only buy, in which case $e_i\leq 0$ and $e_i=a_i-x_i$. 
\end{itemize}

Let $\lambda^\ast$ continue to be the price for a unit of shared resource. Denote $\mathbf{e}=(e_1\dots e_n)^\top$ as the vector representing the traded resource profile across the network. 

\begin{definition}
A triplet of   price-allocation-trade profile ($\lambda^\ast$, $\mathbf{x}^\ast$, $\mathbf{e}^\ast$) is a competitive equilibrium for the MAS-SALTD if the following conditions hold:

(i) Each agent $i$ maximizes her combined payoff at $(\mathbf{x}^\ast$, $\mathbf{e}^\ast)$ while meeting the physical constraint, i.e., $(x_i^\ast$, $e_i^\ast)$ is an optimizer for the following constrained optimization problem:
\begin{equation}
\begin{aligned}
\max_{{x}_i, e_i} \quad &  f_i(x_i)+\lambda^\ast e_i \\
{\rm s.t.} \quad & x_i+e_i\leq a_i \\
\quad & x_i\in \mathbb{R}^{\geq0}, e_i\in \mathbb{R}.
\end{aligned}\label{trading_ceq}
\end{equation}

(ii) The total demand and supply are balanced across the network:
\begin{equation}\label{trading_demand_supply_constraints}
	\sum_{i=1}^n e_i^\ast =0.
\end{equation} 
\end{definition}

\begin{definition}
A pair of resource allocation-trade profile $(\mathbf{x}^\star, \mathbf{e}^\star)$ is a social welfare  equilibrium for the MAS-SALTD if it is an optimizer  to the following optimization problem:  
\begin{align}
\max_{\mathbf{x},\mathbf{e}} \quad & \sum_{i=1}^n f_i(x_i)\label{eq11} \\
{\rm s.t.} \quad & \sum_{i=1}^n e_i = 0,\label{eq12}\\
\quad & x_i+e_i\leq a_i; i\in\mathrm{V},\label{eq13} \\
\quad & x_i\in \mathbb{R}^{\geq0}, e_i\in \mathbb{R}; i\in\mathrm{V}.\label{eq14}
\end{align}

\end{definition}

\begin{theorem}\label{thm2}Consider the MAS-SALTD. 
Suppose each $f_i(\cdot)$ is   concave over the domain $\mathbb{R}^{\geq 0}$. Then the social welfare  equilibrium(s)  and the competitive equilibrium(s)  continue to coincide under the shared load decisions for the agents. To be precise, the following statements hold.

(i) If  ($\lambda^\ast, \mathbf{x}^\ast, \mathbf{e}^\ast$) is a competitive equilibrium, then $(\mathbf{x}^\ast, \mathbf{e}^\ast)$ is a social welfare  equilibrium. 

(ii) If $(\mathbf{x}^\star, \mathbf{e}^\star)$ is a social welfare  equilibrium, then there exists $\lambda^\ast \in \mathbb{R}$ such that ($\lambda^\ast, \mathbf{x}^\star, \mathbf{e}^\star$) is a competitive equilibrium.

\end{theorem}
{\it Proof.} (i) The proof follows the same analysis as the proof of Theorem 1(i), where the desired connection is in place with the definitions of the optimization goals, respectively, for the competitive equilibrium  and the social welfare equilibrium. 

(ii) 
The key idea of the proof continues to be based on strong duality applied to in the definition of the social welfare equilibrium, as the proof of Theorem 1. However, now the social welfare equilibrium contains additional inequality constraints $x_i+e_i\leq a_i, i\in \mathrm{V}$. Inevitably they will lead to auxiliary dual variables, in addition to the dual variable related to the equality constraint $\sum e_i=0$ if we simply repeat the proof of Theorem 1. In order to highlight the role of the dual variable corresponding to the equality constraint, and establish it as the price in competitive equilibrium, we need a refined treatment. To this end, we define a set $\mathbb{X}_i$ for all $i \in \mathrm{V}$ in terms of the inequality constraint in \eqref{eq13} as $\mathbb{X}_i=\{(x_i,e_i)|x_i+e_i \leq a_i; x_i\in \mathbb{R}^{\geq 0};e_i\in \mathbb{R}\}$. Clearly, $\mathbb{X}_i$ is a polyhedral set (see Chapter 3.4.2, Duality Theory in \cite{Bert_N&P}). Denoting again $f_i=-g_i$, the problem  \eqref{eq11}-\eqref{eq14} can be written as:
\begin{equation}\label{eq15}
\begin{aligned}
\max  \quad &  \sum_{i=1}^n g_i(x_i)\\
{\rm s.t.} \quad & (x_i,e_i) \in \mathbb{X}_i, i\in \mathrm{V} \\
\quad &  \sum_{i=1}^n e_i =0.\
\end{aligned} 
\end{equation}

Let $\tau$ be the Lagrange multiplier associated with constraint $\sum_{i=1}^n e_i =0$. Subsequently, we can define the dual function where the primal variables are in a polyhedral set as ( \cite{Bert_N&P}, section 5.1.6): 
\begin{equation}\label{eq17}
L^*(\tau)=\sum_{i=1}^n L^*_i(\tau) ,
\end{equation}
where
\begin{equation}
    L^*_i(\tau)=\inf_{(x_i,e_i) \in \mathbb{X}_i} \Big(g_i(x_i)+\tau e_i\Big), \hspace{0.5cm} i \in \mathrm{V}.
\end{equation}
Let $(\mathbf{x}^\star, \mathbf{e}^\star)$ be a social welfare equilibrium and $\tau^\ast$ be the dual optimal i.e. $\tau^\ast \in \argmin\limits_{\tau \in \mathbb{R}} L^*(\tau)$. Since the problem \eqref{eq15}  is feasible and its optimal value is finite, strong duality holds (\cite{Bert_N&P}, Proposition 5.2.1). This means that 
\begin{align}
    \sum_{i=1}^n g_i(x_i^\star)&= L^*(\tau^\ast)\label{eq19}\\
                               &= \sum_{i=1}^n \left(\inf_{(x_i,e_i) \in \mathbb{X}_i}\Big(g_i(x_i)+\tau^\ast e_i\Big)\right)\label{eq20}\\
                               &\leq \sum_{i=1}^n g_i(x_i^\star)+\tau^\ast \sum_{i=1}^n e_i^\star\label{eq21}\\
                               &\leq \sum_{i=1}^n g_i(x_i^\star)\label{eq22}.
\end{align}
Equation \eqref{eq19} states that the duality gap is zero, \eqref{eq20} comes from the definition of the dual function, \eqref{eq21} follows since the minimization of $\sum_{i=1}^n g_i(x_i)+\tau^\ast \sum_{i=1}^n e_i$ over $(x_i,e_i) \in \mathbb{X}_i$ is always less than or equal to the value at $\sum_{i=1}^n g_i(x_i^\star)+\tau^\ast \sum_{i=1}^n e_i^\star$, \eqref{eq22} follows from $\sum_{i=1}^n e_i^\star=0$. We conclude that the two inequalities hold with equality which implies $(\mathbf{x}^\star, \mathbf{e}^\star)$ minimizes $\sum_{i=1}^n g_i(x_i)+\tau^\ast \sum_{i=1}^n e_i$ over $(x_i,e_i) \in \mathbb{X}_i$. Therefore, there holds
\begin{align}\label{eq23}
    (\mathbf{x}^\star, \mathbf{e}^\star) \in \argmin_{\substack{(x_i,e_i) \in \mathbb{X}_i,\\ i\in \mathrm{V}}} \sum_{i=1}^n g_i(x_i)+\tau^\ast \sum_{i=1}^n e_i.
\end{align}
Since \eqref{eq23} is separable in all $i \in \mathrm{V}$, an equivalent formulation is 
\begin{align}\label{eq24}
    (x_i^\star, e_i^\star) \in \argmin_{(x_i,e_i) \in \mathbb{X}_i} g_i(x_i)+\tau^\ast e_i,\quad i \in \mathrm{V}.
\end{align}
Let us define the equilibrium price $\lambda^\ast$ as $\lambda^\ast=-\tau^\ast$. It follows from \eqref{eq24} that $(x_i^\star, e_i^\star)$ is the solution of the following optimization problem:
\begin{equation} 
\begin{aligned}
\max   \quad &  f_i(x_i)+\lambda^\ast e_i\\
{\rm s.t.} \quad & x_i+e_i \leq a_i \\
\quad &  x_i\in \mathbb{R}^{\geq0}, e_i\in \mathbb{R}.\
\end{aligned} 
\end{equation}
Hence, we conclude that the triplet ($\lambda^\ast, \mathbf{x}^\star, \mathbf{e}^\star$) is a competitive equilibrium.  \hfill$\square$

In the presence of agent trading decisions, the price $\lambda^\ast$ under any competitive equilibrium must be non-negative, as shown in the following example. 
\begin{proposition}\label{prop2} Consider the MAS-SALTD. 
Suppose each $f_i(\cdot)$ is   concave over the domain $\mathbb{R}^{\geq 0}$. Let ($\lambda^\ast, \mathbf{x}^\ast, \mathbf{e}^\ast$) be a competitive equilibrium under the agent trading decisions.  Then there always holds that $\lambda^\ast\geq 0$.  
\end{proposition}
{\it Proof.} Assume $\lambda^\ast<0$. Then $f_i(x_i) + \lambda^\ast e_i$ is a strictly decreasing function with respect to $e_i$. Since $e_i$ is unbounded below and upper bounded by $e_i \leq a_i-x_i$, there can not be a finite ${e}^\ast_i$ such that  ${e}^\ast_i\in \arg \max_{e_i\in \mathbb{R}}\big( f_i(x_i) + \lambda^\ast e_i\big)$, contradicting the definition of the competitive equilibrium. This completes the proof. \hfill$\square$

\section{Social Shaping for Competitive Equilibrium}\label{sec:shaping}
Consistent with classical welfare economics theory, a competitive equilibrium, despite being a social welfare equilibrium as well, indicates nothing about fairness or sustainability. If the optimal pricing $\lambda^\ast$ is too high, agents would option of the system, instead of participating in the self-sustained multi-agent system. When members leave the system, the achievable payoff for the remaining agents would go down.  Therefore, the agents share a social responsibility in shaping their utility functions so that $\lambda^\ast$ is within a socially acceptable range.   

\subsection{Shaping the  Competitive Equilibrium}
Now we present an approach to achieve a socially acceptable competitive equilibrium, by synthesizing a class of utility functions  from which agents can select.  We make the following assumption.  

\medskip

\noindent{\bf Assumption 1.} Each $f_i$ is represented by $f_i(x_i)= -\frac{1}{2}b_ix_i^2 +k_i x_i$, where $b_i\in\mathbb{R}^{\geq 0}$ and $k_i\in\mathbb{R}^{\geq 0}$; a utility function $f_i$ is socially admissible  if   there hold $k_i\in[k_{\rm min},k_{\rm max}]$ and $b_i\in[b_{\rm min},b_{\rm max}]$. 

\medskip

Let $\lambda^\dag >0$ represent the highest pricing for $\lambda^\ast$ that agents can accept, and we term such a competitive equilibrium $\lambda^\ast\leq\lambda^\dag $ a {\em socially resilient equilibrium}. Let $\mathbf{a}=(a_1\dots a_n)^\top$ represent the  network resource allocation profile, and let $C:=\sum_{i=1}^n a_i$ represent the network resource capacity. Assuming $C$ and $\mathbf{a}$ are given network characteristics, we consider the following problem of shaping the  competitive equilibrium. 

\noindent{\bf Problem.} (Social Competitive Equilibrium Shaping)  Consider the MAS-SALD.  Find the range for $k_{\rm min}$, $k_{\rm max}$, $b_{\rm min}$, $b_{\rm max}$ under which  there always exists a competitive equilibrium that leads to $0\leq \lambda^\ast \leq \lambda^\dag$, 
 for all socially admissible  utility functions. 
 
 \subsection{Socially Admissible  Utility Functions}

Denote $\mathbf{k}=(k_1,\dots,k_n)$ and $\mathbf{b}=(b_1,\dots,b_n)^\top$. For two vectors $\mathbf{l}=(l_1,\dots,l_n)$ and $\mathbf{l}'=(l_1',\dots,l_n')$, we write $\mathbf{l}	\preceq  \mathbf{l}'$ if there holds $l_i\leq l_i'$ for all $i=1,\dots,n$. In other words, $\preceq$ defines a partial order for all vectors in $\mathbb{R}^n$.   

Define 
\begin{align} \label{eq25}
\mathscr{S}_\ast:= \Big\{ \big(k_{\rm min}, k_{\rm max}, b_{\rm min}, b_{\rm max}\big)\in \mathbb{R}_{\geq 0}^4: \frac{n k_{\rm min}}{b_{\rm max}} \geq C ; -\frac{n k_{\rm min}}{b_{\rm max}} + \frac{n k_{\rm max}}{b_{\rm min}} \leq C ;  -\frac{n \lambda^{
\dag} }{b_{\rm max}} + \frac{n k_{\rm max}}{b_{\rm min}} \leq C \Big\}.
\end{align}

We present the following theorem.

\begin{theorem}\label{thm3}
Consider the MAS-SALD. Let Assumption 1 hold. The following statements hold. 

(i) The competitive equilibrium is unique, and therefore, there exists a well-defined mapping, denoted by $\mathcal{F}(\cdot,\cdot)$, that maps $(\mathbf{k},\mathbf{b}) $ to  $\lambda^\ast:=\mathcal{F}(\mathbf{k},\mathbf{b})$ where $\lambda^\ast$ belongs to the competitive equilibrium. 

(ii) The competitive equilibrium is always socially resilient (ie $\lambda^\ast\leq\lambda^\dag $) for all socially admissible utility functions as long as $(k_{\rm min}, k_{\rm max}, b_{\rm min}, b_{\rm max}\big)\in \mathscr{S}_\ast$.

(iii) Let $(k_{\rm min}, k_{\rm max}, b_{\rm min}, b_{\rm max}\big)\in \mathscr{S}_\ast$ be given. 
The mapping $ \mathcal{F}(\cdot,\cdot)$ is monotone  under the partial order $\preceq$ over $\mathbf{k}$   in the sense that
$$
\mathcal{F}(\mathbf{k},\mathbf{b}) \leq  \mathcal{F}(\mathbf{k}',\mathbf{b})
$$ 
for all socially admissible $\mathbf{k}'	\preceq  \mathbf{k}$.
 \end{theorem}
{\it Proof.} (i) Under Assumption 1, all $f_i$ are strictly concave. Thus, the optimization  problem (\ref{eq1}) is strictly convex, leading to a unique optimal solution. According to Theorem 1, in a competitive equilibrium $(\lambda^\ast,\mathbf{x}^\ast)$, $\mathbf{x}^\ast$ must be the unique primal optimal solution for (\ref{eq1}). While from the definition of equilibrium, $\lambda^\ast$ is by definition the optimal dual variable for (\ref{eq1}), which is also unique.

(ii) Let $(\lambda^\ast,\mathbf{x}^\ast)$ be a competitive equilibrium. Then based on the definition of competitive equilibrium, there holds 
$
 {x}_i^\ast 
$
is the optimal solution to 
$$
\max_{{x}_i \in\mathbb{R}_{\geq 0}} -\frac{1}{2}b_ix_i^2 +k_i x_i +\lambda^\ast (a_i-x_i), \quad i = 1,2,\cdots,n .
$$
This implies that under the condition $\lambda^\ast\leq k_i$, there holds
\begin{align}
 {x}_i^\ast=\frac{k_i- \lambda^\ast}{b_i}. 
\end{align}

We substitute the above form of  ${x}_i^\ast$ into the constraint $\sum_{i=1}^{n} x_i =C$ to obtain
\begin{align}
\Big(\sum_{i=1}^n \frac{1}{b_i}\Big)\lambda^\ast=\sum_{i=1}^n \frac{k_i}{b_i} -C .
\end{align}
We thus confirm that 
\begin{align}\label{s1}
0\leq \lambda^\ast &=  \Big(\sum_{i=1}^n \frac{k_i}{b_i} - C\Big)/\Big(\sum_{i=1}^n \frac{1}{b_i}\Big) \nonumber\\
&\leq \Big(  \frac{n k_{\rm max}}{b_{\rm min}} - C\Big)/\Big(\frac{n}{b_{\rm max}}\Big)
\end{align}
if there holds 
$$
\frac{n k_{\rm min}}{b_{\rm max}}  \leq C.
$$
From (\ref{s1}) we also know
\begin{align}
-\frac{n k_{\rm min}}{b_{\rm max}} + \frac{n k_{\rm max}}{b_{\rm min}} \leq C 
\end{align}
guarantees $\lambda^\ast\leq k_i$. Collecting all conditions for $(k_{\rm min}, k_{\rm max}, b_{\rm min}, b_{\rm max}\big)$, we obtain that      $\lambda^\ast\leq\lambda^\dag $  for all socially admissible utility functions with $(k_{\rm min}, k_{\rm max}, b_{\rm min}, b_{\rm max}\big)\in \mathscr{S}_\ast$. 

(iii) We have established that if $(k_{\rm min}, k_{\rm max}, b_{\rm min}, b_{\rm max}\big)\in \mathscr{S}_\ast$, then
\begin{align}
\lambda^\ast =\mathcal{F}(\mathbf{k},\mathbf{b})=  \Big(\sum_{i=1}^n \frac{k_i}{b_i}-C\Big)/\Big(\sum_{i=1}^n \frac{1}{b_i}\Big).
\end{align}
It is straightforward to verify that $$
\mathcal{F}(\mathbf{k},\mathbf{b}) \leq  \mathcal{F}(\mathbf{k}',\mathbf{b})
$$ 
for all socially admissible $\mathbf{k}	\preceq  \mathbf{k}'$. \hfill$\square$

\section{Dynamic  Multi-Agent Systems}\label{sec:dynamic}

\subsection{MAS with Dynamic Agent Load/Trading Decisions}
Here we consider the load balancing problem for dynamical multi-agent systems. 

\medskip

\noindent{\bf MAS with Dynamic Agent Load/Trading Decisions} (MAS-DALTD). Each agent $i\in\mathrm{V}$  is associated with a dynamical state $\mathbf{y}_i(t)\in \mathbb{R}^m$, described by 
\begin{align}
\mathbf{y}_i(t+1)= \mathbf{A}_i \mathbf{y}_i(t)+\mathbf{B}_i \mathbf{u}_i(t), \quad t=0,\dots, T-1,
\end{align}
where $\mathbf{u}_i(t)\in \mathbb{R}^m$ is the control input, and $\mathbf{A}_i$ and $\mathbf{B}_i$ are real matrices with proper dimensions. Associated with $t=0,\dots,T-1$, agent $i$ incurs a utility function $f_i(\mathbf{y}_i(t),\mathbf{u}_i(t))$; the terminal utility for agent $i$ is $\Phi_i(\mathbf{y}_i(T))$. Upon taking the control action $\mathbf{u}_i(t)$, the required resource is defined by the function $h_i(\mathbf{u}_i(t))$. Each agent can produce an $a_i(t)$ units of energy at time $t$, and also makes a trading decision $e_i(t)$ units of energy over  the network at time $t$.  Similarly, \begin{itemize}
\item[(i)] if $h_i(\mathbf{u}_i(t))<a_i(t)$, then agent $i$ can sell, in which case $e_i(t) \geq 0$ and $e_i(t)\leq a_i(t)-h_i(\mathbf{u}_i(t))$;
\item[(ii)] if $h_i(\mathbf{u}_i(t))\geq a_i(t)$, then agent $i$ will buy, in which case $e_i(t)\leq 0$ and $e_i(t)=a_i(t)-h_i(\mathbf{u}_i(t))$. 
\end{itemize}
We denote $\bm{\lambda}=(\lambda_0\dots \lambda_{T-1})^\top$ as the pricing vector through the time horizon, where $\lambda_t$ is the unit price for traded energy at step $t$. Consequently, the payoff of agent $i$ throughout $[0,T]$ is described by 
$$
\sum_{t=0}^{T-1}\bigg( f_i(\mathbf{y}_i(t),\mathbf{u}_i(t))+ \lambda_t e_i(t)\Big)+\Phi(\mathbf{y}_i(T)). 
$$ 

Denote $\mathbf{y}(t)=(\mathbf{y}_1(t)^\top\dots \mathbf{y}_n(t)^\top)^\top$, $\mathbf{u}(t)=(\mathbf{u}_1(t)^\top\dots \mathbf{u}_n(t)^\top)^\top$, and $\mathbf{e}(t)=(\mathbf{e}_1(t)^\top\dots \mathbf{e}_n(t)^\top)^\top$. Further define $\mathbf{Y}=(\mathbf{y}(0)^\top\dots \mathbf{y}(T)^\top)^\top$, $\mathbf{U}=(\mathbf{u}(0)^\top\dots \mathbf{u}(T-1)^\top)^\top$ and $\mathbf{E}=(\mathbf{e}(0)^\top\dots \mathbf{e}(T-1)^\top)^\top$. Also introduce $\mathbf{U}_i=(\mathbf{u}_i(0)^\top\dots \mathbf{u}_i(T-1)^\top)^\top$, $\mathbf{E}_i=(\mathbf{e}_i(0)^\top\dots \mathbf{e}_i(T-1)^\top)^\top$ and $\mathbf{a}_i=(a_i(0)^\top\dots a_i(T-1)^\top)^\top$

\begin{definition}
Let $\mathbf{y}(0)=\mathbf{y}_0 \in\mathbb{R}^{mn}$ be given. A  triple of price-control-trading profiles   ($\bm{\lambda}^\ast,\mathbf{U}^\ast,\mathbf{E}^\ast$) is a dynamic competitive equilibrium if the following conditions hold:

(i) each agent $i$ maximizes its combined payoff under $\mathbf{U}_i^\ast$ and $\mathbf{E}^\ast_i$:
\begin{equation}
\begin{aligned}
\max_{\mathbf{U}_i,\mathbf{E}_i} \quad & \sum_{t=0}^{T-1}\bigg( f_i(\mathbf{y}_i(t),\mathbf{u}_i(t))+ \lambda_t^\ast e_i(t)\Big)+\Phi(\mathbf{y}_i(T))  \\
{\rm s.t.} \quad & \mathbf{y}_i(t+1)= \mathbf{A}_i \mathbf{y}_i(t)+\mathbf{B}_i \mathbf{u}_i(t),\\
           \quad &  e_i(t)\leq a_i(t)-h_i(\mathbf{u}_i(t)),t=0,\dots,T-1;
\end{aligned} \label{eq102}
\end{equation}

(ii) the total demand and supply are balanced across the network for all time, i.e., there holds $$
\sum_{i=1}^n e_i(t)=0, \quad t=0,\dots,T-1.
$$  
\end{definition}

\begin{definition}
Let $\mathbf{y}(0)=\mathbf{y}_0 \in\mathbb{R}^{mn}$ be given. A pair of control-trading  profiles $(\mathbf{U}^\star,\mathbf{E}^\star)$ is a dynamic social welfare  equilibrium if it is a solution to the following optimal control problem:  
\begin{align}
\max_{\mathbf{U},\mathbf{E}} \quad & \sum_{i=1}^n \Big(\sum_{t=0}^{T-1} f_i(\mathbf{y}_i(t),\mathbf{u}_i(t))  
+\Phi(\mathbf{y}_i(T)) \Big) \label{dym_obj} \\
{\rm s.t.} \quad &  \mathbf{y}_i(t+1)= \mathbf{A}_i \mathbf{y}_i(t)+\mathbf{B}_i \mathbf{u}_i(t), \quad t=0,\dots, T-1, \; i\in\mathrm{V}  \label{dym_con1} \\
 \quad &  e_i(t)\leq a_i(t)-h_i(\mathbf{u}_i(t)),  \quad t=0,\dots, T-1, \; i\in\mathrm{V}  \label{dym_con2}\\
\quad & \sum_{i=1}^n e_i(t) =0, \quad t=0,\dots,T-1. \label{dym_con3}
\end{align}
\end{definition}

\subsection{Dynamic Competitive Equilibrium}
We impose the following assumption.

\medskip

\noindent{\bf Assumption 2}. (i) $\Phi$ is a concave function; (ii) the $f_i$  are concave functions  for  $i\in\mathrm{V}$; (iii) the $h_i$  are nonnegative convex functions  for  $i\in\mathrm{V}$, and $h_i(\mathbf{z})<b$ defines a bounded open set of $\mathbf{z}$ in $\mathbb{R}^m$ for $b>0$; (iv) $\sum_{i=1}^n a_i(t) >0$ for all $t=1,\dots,T-1$. 

\medskip

We present the following result which establishes the similar connection between the competitive equilibrium and social welfare equilibrium under this dynamic setting. 

\begin{theorem}\label{thm4} Consider the  MAS-DALTD with $\mathbf{y}(0)=\mathbf{y}_0 \in\mathbb{R}^{mn}$ be given. Let Assumption 2 hold. The dynamic social welfare  equilibrium(s)  and the dynamic competitive equilibrium(s) coincide and the following statements hold.

(i) If  $(\bm{\lambda}^\ast,\mathbf{U}^\ast,\mathbf{E}^\ast)$ is a dynamic competitive equilibrium, then $(\mathbf{U}^\ast,\mathbf{E}^\ast)$ is a dynamic social welfare  equilibrium. 

(ii) If $(\mathbf{U}^\star,\mathbf{E}^\star)$ is a dynamic social welfare  equilibrium, then there exists $\bm{\lambda}^\ast \in \mathbb{R}^T$ such that $(\bm{\lambda}^\ast,\mathbf{U}^\star,\mathbf{E}^\star)$ is a competitive equilibrium.

\end{theorem}
\textit{Proof:} (i) The proof of sufficiency follows from a similar analysis as the proof of Theorem \ref{thm1}. The transition from competitive equilibrium to social welfare equilibrium under this dynamical setting continues to be a direct consequence of the formulations of the two underlying optimization problems.

(ii) First of all, the dynamics $\mathbf{y}_i(t+1)= \mathbf{A}_i \mathbf{y}_i(t)+\mathbf{B}_i \mathbf{u}_i(t)$ with given $\mathbf{y}(0)$ ensures that any $\mathbf{y}_i(t)$ for $t=1,\dots,T$ is a linear combination of $\mathbf{y}(0)$ and $\mathbf{u}_i(0),\dots, \mathbf{u}_i(t-1)$. Therefore, we can always write for any $i\in\mathrm{V}$ that
\begin{align}
\mathbf{y}_i(t)= p_{i,t}\mathbf{y}_0 +q_{i,t}\mathbf{U}_i,\ t=0,\dots,T
\end{align}
with $p_{i,t}$ and $q_{i,t}$ being matrices with proper dimensions. As a result, we have
\begin{align}
f_i(\mathbf{y}_i(t),\mathbf{u}_i(t)) =f_i(p_t\mathbf{y}_0 +q_t\mathbf{U}_i,\mathbf{u}_i(t)):=\tilde{f}_{i,t}(\mathbf{U}_i). 
\end{align}
In view of Assumption 2, and the fact that composition of a concave function and an affine function continues to be concave, we conclude that  $g_{i,t}(\cdot)$ is a concave function. Similarly, 
$$
\Phi(\mathbf{y}_i(t))=\Phi(p_{i,T}\mathbf{y}_0 +q_{i,T}\mathbf{U}_i):=\Phi_{i} (\mathbf{U}_i) 
$$
where $\Phi_{i}(\cdot)$ is a concave function. 

The optimization problem (\ref{dym_obj})-(\ref{dym_con3}) can be equivalently rewritten as the following convex programming:
\begin{equation}
\begin{aligned}
\min_{\mathbf{U},\mathbf{E}} \quad &  -\sum_{i=1}^n \Big(\sum_{t=0}^{T-1} \tilde{f}_{i,t}(\mathbf{U}_i) 
+\Phi_{i} (\mathbf{U}_i) \Big)\\
{\rm s.t.} \quad   & h_i(\mathbf{u}_i(t)) +e_i(t)\leq a_i(t),  \quad t=0,\dots, T-1, \; i\in\mathrm{V} \\
\quad & \sum_{i=1}^n e_i(t) =0, \quad t=0,\dots,T-1.
\end{aligned} \label{103}
\end{equation}
Similarly, (\ref{eq102}) can be equivalently written as convex programming 
\begin{equation}
\begin{aligned}
\min_{\mathbf{U}_i,\mathbf{E}_i} \quad & - \sum_{t=0}^{T-1}\Big( \tilde{f}_{i,t}(\mathbf{U}_i)+ \lambda_t^\ast e_i(t)\Big)
+\Phi_{i} (\mathbf{U}_i)  \\
{\rm s.t.} \quad &  e_i(t)\leq a_i(t)-h_i(\mathbf{u}_i(t)),t=0,\dots,T-1.
\end{aligned} \label{eq105}
\end{equation}

Next, with  Assumption 2.(iii)-(iv), we can verify that Slater's condition holds for (\ref{103}) and (\ref{eq105}),  which guarantees strong duality in both problems \cite{boyd}. Also, noting 
\begin{align}
\sum_{i=1}^n h_i(\mathbf{u}_i(t)) \leq \sum_{i=1}^n a_i(t)
\end{align}
and the assumption that  $h_i(\mathbf{z})<b$ defines a bounded open set of $\mathbf{z}$ in $\mathbb{R}^m$ for $b>0$, $\mathbf{u}_i(t)$ takes values in a compact set for all $i$ and $t$. Moreover, since $h_i(\mathbf{u}_i(t))\geq 0$ holds for all $i$ and for all $t$, there holds $e_i(t)\leq a_i(t)$. The constraint  $\sum_{i=1}^n e_i(t)=0$ further ensures $e_i(t)\geq -\sum_{i=1}^n a_i(t)$ for all $i$ and $t$. Thus $e_i(t)$ also takes values in a compact set for all $i$ and $t$. The convex programming problem (\ref{103}) leads to finite primal  solution.  
 
The Lagrangian dual function of (\ref{103}) can be written as
\begin{align}\label{106}
L(\mathbf{U},\mathbf{E},\bm{\lambda}, \bm{\mu})&= -\sum_{i=1}^n \Big(\sum_{t=0}^{T-1} \tilde{f}_{i,t}(\mathbf{U}_i) +\Phi_{i} (\mathbf{U}_i) \Big) + \sum_{t=0}^{T-1} \sum_{i=1}^n \lambda_t e_i(t) +\sum_{t=0}^{T-1}\sum_{i=1}^n \mu_{i,t}\big( h_i(\mathbf{u}_i(t)) +e_i(t) - a_i(t) \big)\nonumber\\
&=\sum_{i=1}^n  L_i(\mathbf{U}_i,\mathbf{E}_i,\bm{\lambda}, \bm{\mu}_i)  
\end{align} 
where 
\begin{align}
 L_i(\mathbf{U}_i,\mathbf{E}_i,\bm{\lambda}, \bm{\mu}_i)=-\sum_{t=0}^{T-1} \tilde{f}_{i,t}(\mathbf{U}_i) +\Phi_{i} (\mathbf{U}_i)+\sum_{t=0}^{T-1}   \lambda_t e_i(t)+\sum_{t=0}^{T-1}  \mu_{i,t}\big( h_i(\mathbf{u}_i(t)) +e_i(t) - a_i(t) \big).
\end{align}
Here $\mu_{i,t}\geq 0$ since they correspond to the inequality constraints. We have used the conventional notation $\bm{\mu}_i=(\mu_{i,0},\dots, \mu_{i,T-1})^\top$ and $\bm{\mu}=(\bm{\mu}_1^\top,\dots, \bm{\mu}_n^\top)^\top$. 

Finally, letting an optimal dual solution of (\ref{103}) be $(\bm{\lambda}^\ast, \bm{\mu}^\ast)$, there holds from strong duality 
\begin{align}
 (\mathbf{U}^\star,\mathbf{E}^\star)\in \arg \min\   L(\mathbf{U},\mathbf{E},\bm{\lambda}^\ast, \bm{\mu}^\ast)
\end{align}
if $(\mathbf{U}^\star,\mathbf{E}^\star)$ is a dynamic social welfare  equilibrium. This implies from (\ref{106}) that  
\begin{align}\label{108}
 (\mathbf{U}_i^\star,\mathbf{E}_i^\star)\in \arg \min\   L_i(\mathbf{U},\mathbf{E},\bm{\lambda}^\ast, \bm{\mu}^\ast).
\end{align}
Now, $\bm{\mu}^\ast$ is obtained by solving 
\begin{align}
\max_{\bm{\lambda}, \bm{\mu}}\  \min_{\mathbf{U},\mathbf{E}}\  L(\mathbf{U},\mathbf{E},\bm{\lambda} , \bm{\mu}) =\max_{\bm{\lambda}, \bm{\mu}}\  \min_{\mathbf{U},\mathbf{E}}\ \sum_{i=1}^n  L_i(\mathbf{U}_i,\mathbf{E}_i,\bm{\lambda}, \bm{\mu}_i)  
\end{align}
where the maximization and minimization are taken in their respective domains for $\bm{\lambda},  \bm{\mu}, \mathbf{U}, \mathbf{E}$. As a result, there must hold
\begin{align}\label{109}
\bm{\mu}_i^\ast \in \arg \max_{\bm{\mu}_i} \min_{\mathbf{U}_i,\mathbf{E}_i} L_i(\mathbf{U}_i,\mathbf{E}_i,\bm{\lambda}^\ast, \bm{\mu}_i).
\end{align}
It is worth emphasizing that $L_i(\mathbf{U}_i,\mathbf{E}_i,\bm{\lambda}^\ast, \bm{\mu}_i)$ is, precisely, the Lagrangian of (\ref{eq105}).  Therefore, (\ref{109}) ensures that  $\bm{\mu}_i^\ast$ is an optimal dual solution of (\ref{eq105}),  and then from strong duality (\ref{108}) ensures that $(\mathbf{U}^\star,\mathbf{E}^\star)$ is an optimal primal solution of (\ref{eq105}). In other words, we have proven $(\bm{\lambda}^\ast,\mathbf{U}^\star,\mathbf{E}^\star)$ is a competitive equilibrium.

The proof of the desired theorem is now complete. \hfill$\square$

\section{Numerical Examples}\label{sec:numerical}

\subsection{MAS with Static Decisions}
In this section, several numerical examples are provided for validation of Theorem \ref{thm1}, \ref{thm2}, \ref{thm3} and Proposition \ref{prop2}. These experiments are implemented using the Python software package CVXPY.

\medskip

\noindent{\bf Example 1.} Consider a multi-agent system with four agents who have local resource $(a_1,a_2,a_3,a_4)= (13,14,4,7)$. Each agent $i$  is associated with a utility function $f_{i}$ which is represented by $f_i(x_i)= \min(k_{i}x_{i}, \beta_{i}) $ with $(k_1,k_2,k_3,k_4) = (21,20,23,32)$  and  $(\beta_1,\beta_2,\beta_3,\beta_4) = (135,600,130,150)$.  

(i) Let the multi-agent system be  MAS-SALD. The social welfare equilibrium can be computed by numerically  solving the optimization problem \eqref{eq1} as $$\mathbf{x}^{\star} = (6.429, 21.232, 5.652, 4.688)^{\top},$$ and the corresponding optimal  dual variable is also obtained as $\lambda^{\ast} = 20$. Letting $\lambda^{\ast}=20$, we then compute a competitive equilibrium that satisfies  \eqref{load_competitive_eq}-\eqref{load_demand_supply_constraints} as  $$\mathbf{x}^{\ast} = (6.429, 21.232, 5.652, 4.688)^{\top}.$$ In particular, we obtain $x_1^\ast=6.429$, $x_3^\ast=5.652$, and $x_4^\ast= 4.688$ by solving \eqref{load_competitive_eq}, and further establish $x_2^\ast=21.232$ from \eqref{load_demand_supply_constraints}.  Clearly there holds $\mathbf{x}^{\star} = \mathbf{x}^{\ast}$, which is consistent with Theorem \ref{thm1} .

(ii) Let the multi-agent system be MAS-SALTD. We compute the social welfare equilibrium $(\mathbf{x}^{\star},\mathbf{e}^{\star})$ by    solving the optimization problem  \eqref{eq11}-\eqref{eq14} as $$\mathbf{x}^{\star} = (6.429, 21.232, 5.652, 4.688)^{\top} \quad \mathbf{e}^{\star} = (6.571, -7.232, -1.652, 2.313)^{\top}.$$ The optimal dual variable $\tau^{\ast}$ corresponding to  the equity constraint \eqref{eq12} can be obtained as $\tau^{\ast} = -20$. We take $\lambda^{\ast} = - \tau^{\ast} =20 $ and establish a competitive equilibrium that satisfies \eqref{trading_ceq}-\eqref{trading_demand_supply_constraints} as $$\mathbf{x}^{\ast} = (6.429, 21.232, 5.652, 4.688)^{\top} \quad \mathbf{e}^{\ast} = (6.571, -7.232, -1.652, 2.313)^{\top}.$$ In particular, we compute $(x_{1}^{\ast},e_{1}^{\ast}) = (6.429,6.571)$, $(x_{3}^{\ast},e_{3}^{\ast}) = (5.652,-1.652)$ and $(x_{4}^{\ast},e_{4}^{\ast}) = (4.688,2.313)$  by solving \eqref{trading_ceq}, and further obtain $(x_{2}^{\ast},e_{2}^{\ast}) = (21.232,-7.232)$ from \eqref{trading_demand_supply_constraints}.  Again there holds $(\mathbf{x}^{\star}, \mathbf{e}^{\star}) = (\mathbf{x}^{\ast}, \mathbf{e}^{\ast})$, which validates Theorem \ref{thm2}. \hfill$\square$

\medskip

\noindent{\bf Example 2.}
Consider a multi-agent system with four agents. The utility function for agent $i$ is in the quadratic form $f_{i} = -\frac{1}{2}b_{i}x_{i}^{2}+k_{i}x_{i}$ for $i=1,2,3,4$. We  consider two pairs of system parameters
\begin{equation}
    \mathbf{b} = (2,5,3,4)^{\top} \quad \mathbf{k} = (21,17,23,13)^{\top}; \tag{PM.1}
\end{equation}
\begin{equation}
    \mathbf{b}^{'} = (2,5,3,4)^{\top} \quad \mathbf{k}^{'} = (25,22,24,14)^{\top}\tag{PM.2}.
\end{equation}
Let the network resource capacity   $C=\sum_{i=1}^4 a_i$ take values in an interval $(0,40).$ We sample the interval $(0,40)$ uniformly with a step-size $0.8$ to obtain $50$ different values for $C$. For each $C$,  we compute the optimal prices of the system under MAS-SALD and MAS-SALTD. 

\begin{figure}[htbp]
    \centering
    \includegraphics[width=0.5\textwidth]{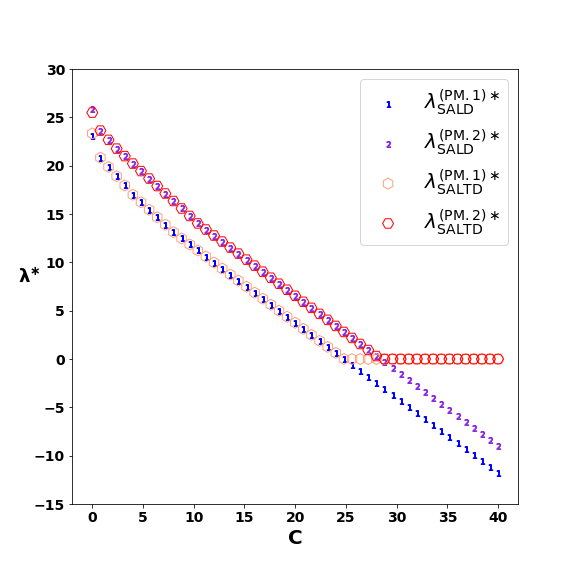}
    \caption{The curves of the optimal prices as functions of the network resource capacity in Example 2. The optimal prices of the MAS-SALD with (PM.1) and (PM.2) are represented by $\lambda^{\rm (PM.1) \ast }_{\rm SALD},\lambda^{\rm (PM.2) \ast }_{\rm SALD}$ respectively, while $\lambda^{\ast \rm (PM.1)}_{\rm SALTD}, \lambda^{\ast \rm (PM.2)}_{\rm SALTD}$ are the optimal prices of MAS-SALTD with (PM.1) and (PM.2) respectively. }
    \label{fig:optimal prices}
\end{figure}

For MAS-SALD, the optimal  dual variables $\lambda^{\ast \rm (PM.1)}_{\rm SALD}$ and $\lambda^{\ast \rm (PM.2)}_{\rm SALD}$ are computed for $10^4$ times corresponding to each value of $C$ by solving \eqref{eq1}, respectively, under the parameter setting  (PM.1) and (PM.2).  For MAS-SALTD,  the optimal  dual variables $\tau^{\ast \rm (PM.1)}_{\rm SALTD}$ and $\tau^{\ast \rm (PM.2)}_{\rm SALTD}$ related  to the equity constraint \eqref{eq12} are also computed for $10^4$ times corresponding to each value of $C$ by solving  \eqref{eq11}-\eqref{eq14}, respectively, under the parameter setting  (PM.1) and (PM.2), and then we take $\lambda^{\ast \rm (PM.1)}_{\rm SALTD} = -\tau^{\ast \rm (PM.1)}_{\rm SALTD}$ and $\lambda^{\ast \rm (PM.2)}_{\rm SALTD} = -\tau^{\ast \rm (PM.2)}_{\rm SALTD}$.
In Figure~\ref{fig:optimal prices}, we plot the $50$ points of optimal prices versus  $C$, to obtain an approximate trajectory of the optimal price as a function of $C$.

From the plot, we easily observe that the optimal price $\lambda^{\ast}_{\rm SALD}$ in MAS-SALD  can indeed take negative values; while the optimal price $\lambda^{\ast}_{\rm SALTD}$  in MAS-SALTD is always non-negative. These  observations are consistent with Proposition \ref{prop1} and Proposition \ref{prop2}. Moreover, for both MAS-SALD and MAS-SALTD, the plot shows that the optimal prices $\lambda^{\ast}_{\rm SALD}, \lambda^{\ast}_{\rm SALTD}$ are decreasing as the network resource capacity $C$ increases. \hfill$\square$
   
\medskip

\noindent{\bf Example 3.}
Consider a MAS-SALD with three agents and    network  capacity $C = 18$. Each agent's utility function is set as the quadratic form $f_{i} = -\frac{1}{2}b_{i}x_{i}^{2}+k_{i}x_{i}$ for $i=1,2,3$. The system's highest pricing for $\lambda^{\ast}$ that agents can accept socially  is assumed to be $\lambda^{\dag} = 42.$ Take $b_{\rm \min}=4$, $b_{\rm \max}=6$, $k_{\rm \min}=40$, and $k_{\rm \max}=50$. We can verify such a configuration of $(b_{\rm \min}, b_{\rm \max}, k_{\rm \min}, k_{\rm \max})$ is a point in  $\mathscr{S}_\ast$ defined in \eqref{eq25}.  

\begin{figure}[htbp]
     \centering
     \begin{subfigure}{0.49\textwidth}
         \centering
         \includegraphics[width=1.2\textwidth]{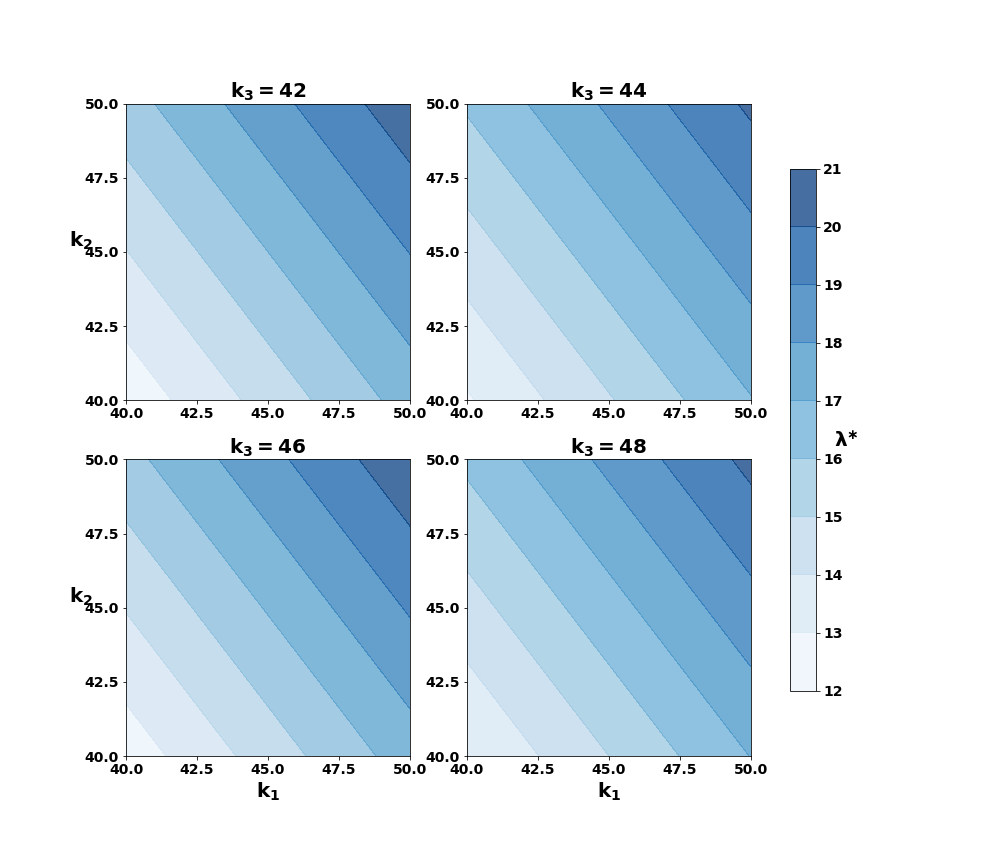}
         \caption{}
         \label{k contour}
     \end{subfigure}
     \begin{subfigure}{0.49\textwidth}
         \centering
         \includegraphics[width=1.2\textwidth]{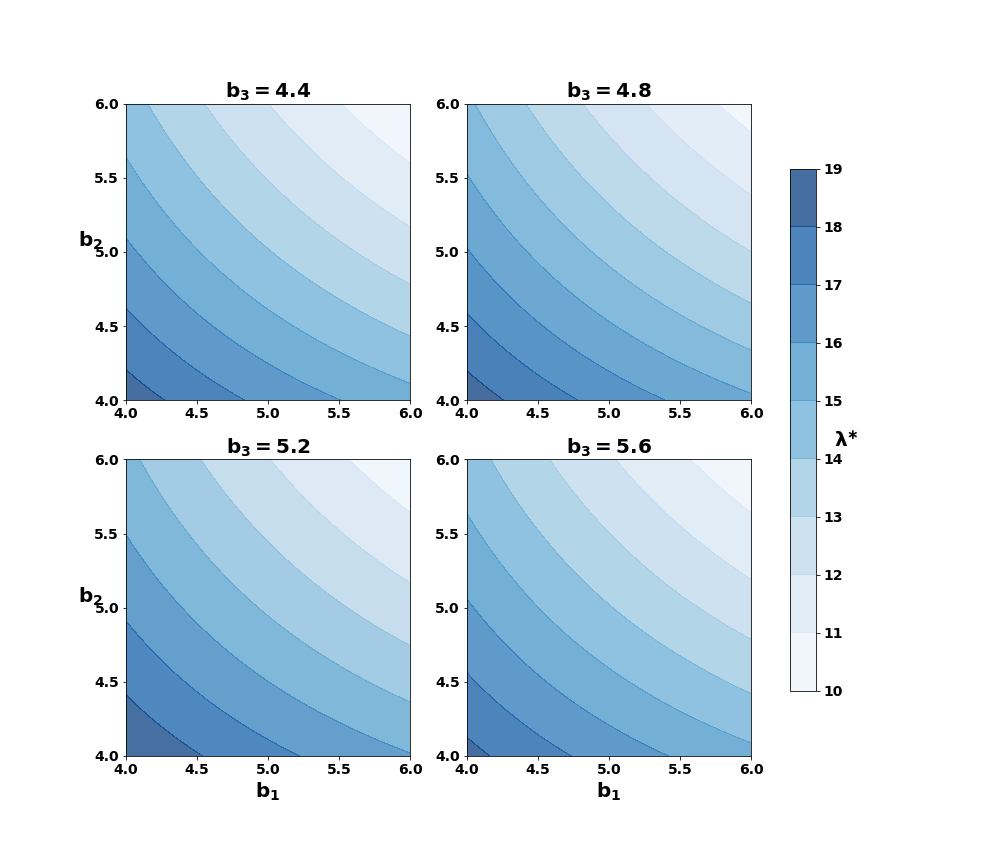}
         \caption{}
         \label{b contour}
     \end{subfigure}
        \caption{Contour maps for the optimal prices in Example 3.}
        \label{Contours}
\end{figure}

(i)  Let $\mathbf{b}$ be fixed to be $\mathbf{b} = (4,5,6)^{\top}$. Take $k_3=42, 44, 46, 48$, respectively.  We sample the space for  $(k_1,k_2)\in[40,50]^2$ and compute the optimal pricing  $\lambda^{\ast}$   by solving the optimal dual variable of \eqref{eq1}. Then we plot the contour maps for the optimal price as a function of $k_1$ and $k_2$ in Figure~\ref{k contour}. 

(ii) Let $\mathbf{k}$ be fixed to be $\mathbf{k} = (44,46,48)^{\top}$. Take $b_3=4.4, 4.8, 5.2, 5.6$, respectively.  We sample the space for  $(b_1,b_2)\in[4,6]^2$ and compute the optimal pricing  $\lambda^{\ast}$   by solving the optimal dual variable of \eqref{eq1}. Then we plot the contour maps for the optimal price as a function of $b_1$ and $b_2$ in Figure~\ref{b contour}. 

From these plots, we observe  that the maximum   value for the price $\lambda^{\ast}$ is $21$, which is lower than $\lambda^{\dag} = 42$. This illustrates all socially admissible utility functions for parameters in the set  $\mathscr{S}_\ast$ lead to socially acceptable prices, providing a validation for Theorem \ref{thm3}(ii). From   Figure~\ref{k contour}, the optimal price is monotone under the partial order $\preceq$  with respect to $\mathbf{k}$, which is consistent with Theorem \ref{thm3}(iii). \hfill$\square$

\subsection{MAS with Dynamic  Decisions}

\begin{figure}[htbp]
    \centering
    \includegraphics[width=0.9\textwidth]{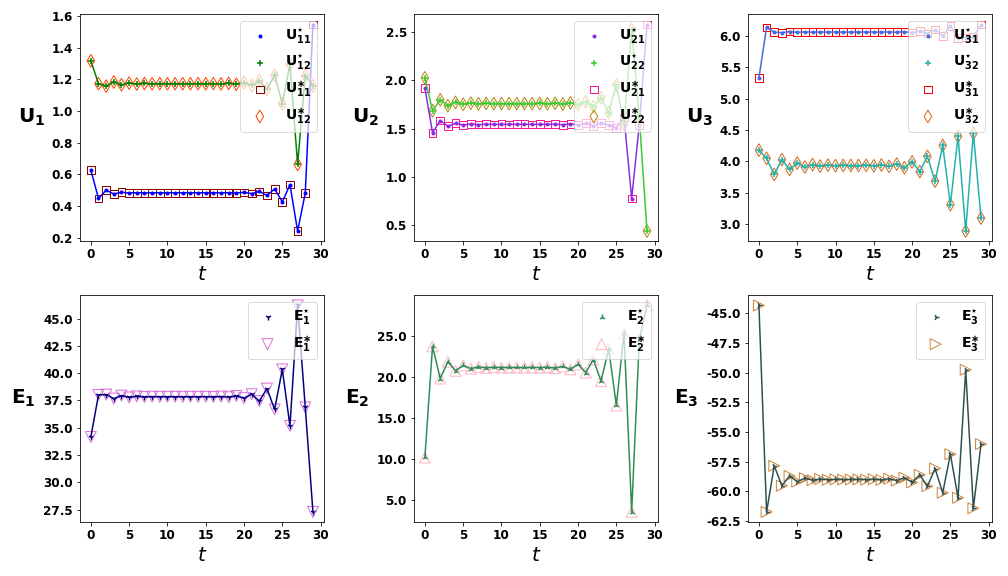}
    \caption{The dynamic social welfare equilibrium and competitive equilibrium in Example 4.}
    \label{fig:dym_swe_ce}
\end{figure}

\begin{figure}
    \centering
    \includegraphics[width=0.5\textwidth]{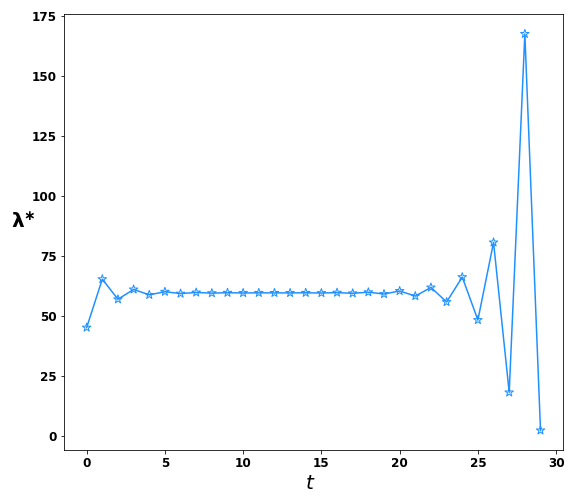}
    \caption{The dynamic optimal price for traded resource versus $t$ in Example 4.}
    \label{fig:dym_price}
\end{figure}

\noindent{\bf Example 4.} Consider a MAS-DALTD with three agents who have initial states $$\mathbf{y}_{1}(0) =\begin{bmatrix}
	 1 \\
	 4 \\
\end{bmatrix}, \mathbf{y}_{2}(0) =\begin{bmatrix}
	 2 \\
	 5 \\\end{bmatrix} ,\mathbf{y}_{3}(0) =\begin{bmatrix}
	 3 \\
	 3 \\\end{bmatrix},$$ and local resource 
	 $$\mathbf{a}_{1} = [50,\cdots,50]_{1\times T}, \quad \mathbf{a}_{2} = [50,\cdots,50]_{1\times T} , \quad
	 \mathbf{a}_{3} = [30,\cdots,30]_{1\times T}.$$ 
The dynamical state $\mathbf{y}_{i}(t) \in \mathbb{R}^{2}$ of agent $i$ is described by $$ \mathbf{y}_i(t+1)= \mathbf{A}_i \mathbf{y}_i(t)+\mathbf{B}_i \mathbf{u}_i(t), \quad t=0,\dots, T-1, $$
where $$\mathbf{A}_1 = 
\begin{bmatrix}
	-0.6 & 0 \\
	 0 & -0.7 \\
\end{bmatrix} , 
\mathbf{A}_2 = \begin{bmatrix}
	-0.5 & 0 \\
	 0 & -0.2 \\
\end{bmatrix} , 
\mathbf{A}_3 = \begin{bmatrix}
	-0.4 & 0 \\
	 0 & -0.8 \\
\end{bmatrix} ,  \mathbf{B}_1 = 
\begin{bmatrix}
	2 & 0\\
	0 & 7\\
\end{bmatrix}, 
\mathbf{B}_2 = \begin{bmatrix}
	4 & 0\\
	0 & 6\\
\end{bmatrix}, 
\mathbf{B}_3 = \begin{bmatrix}
	9 & 0\\
	0 & 3\\
\end{bmatrix}.$$ The utility function of agent $i$ is in the quadratic form $$ f_i(\mathbf{y}_i(t),\mathbf{u}_i(t)) = \mathbf{y}_{i}^{\top}(t) \mathbf{R}_{i} \mathbf{y}_{i}(t) +  \mathbf{W}_{i} \mathbf{y}_{i}(t) + \mathbf{u}_{i}^{\top}(t) \mathbf{Q}_{i} \mathbf{u}_{i}(t) + \mathbf{K}_{i} \mathbf{u}_{i}(t),$$
where $$\mathbf{R}_{1} = 
\begin{bmatrix}
	-5 & 0\\
	0 & -8\\
\end{bmatrix} , 
\mathbf{R}_{2} = 
\begin{bmatrix}
	-3 & 0\\
	0 & -7\\
\end{bmatrix} , 
\mathbf{R}_{3} = 
\begin{bmatrix}
	-2 & 0\\
	0 & -1\\
\end{bmatrix} ,
\mathbf{Q}_{1} = \begin{bmatrix}
	-5 & 0\\
	0  & -4 \\
\end{bmatrix}, 
\mathbf{Q}_{2} = \begin{bmatrix}
	-1 & 0\\
	 0 & -6\\
\end{bmatrix}, 
\mathbf{Q}_{3} = \begin{bmatrix}
	-3 & 0\\
	0 & -2 \\
\end{bmatrix}, $$ $$
\mathbf{W}_1 = \begin{bmatrix}
	200 & 300\\
\end{bmatrix},
\mathbf{W}_2 = \begin{bmatrix}
	200 & 400\\
\end{bmatrix},
\mathbf{W}_3 = \begin{bmatrix}
	450 & 300\\
\end{bmatrix},
\mathbf{K}_1 = \begin{bmatrix}
	50 & 60 \\
\end{bmatrix},
\mathbf{K}_2 = \begin{bmatrix}
	50 & 20\\
\end{bmatrix},
\mathbf{K}_3 = \begin{bmatrix}
	80 & 20\\
\end{bmatrix}.$$ The terminal utility of agent $i$ is also set as the quadratic form 
$$\Phi_{i}(\mathbf{y}_{i}(T)) = \mathbf{y}_{i}^{\top}(T) \mathbf{R}_{i} \mathbf{y}_{i}(T) +  \mathbf{W}_{i} \mathbf{y}_{i}(T). $$
Upon taking $\mathbf{u}_{i}(t),$ the required resource is determined by 
$$h_{i}(\mathbf{u}_{i}(t)) = \mathbf{u}_{i}^{\top}(t) \mathbf{H}_{i} \mathbf{u}_{i}(t), $$
where $$\mathbf{H}_{1} = 
\begin{bmatrix}
	5 & 0\\
	0 & 8\\
\end{bmatrix} ,
\mathbf{H}_{2} = 
\begin{bmatrix}
	3 & 0\\
	0 & 7\\
\end{bmatrix} , 
\mathbf{H}_{3} = 
\begin{bmatrix}
	2 & 0\\
	0 & 1\\
\end{bmatrix}. $$
Let the time horizon take the value of $T = 30.$ We compute the dynamic social welfare equilibrium $(\mathbf{U}^{\star}, \mathbf{E}^{\star})$  by solving the optimization problem (\ref{dym_obj})-(\ref{dym_con3}) and the optimal dual variables $-\bm{\lambda}^{\ast}$ corresponding to \eqref{dym_con3}. 
Given $\bm{\lambda}^{\ast}$, we further compute the dynamic competitive equilibrium $(\mathbf{U}^{\ast}, \mathbf{E}^{\ast})$ by solving \eqref{eq102}. In Figure~\ref{fig:dym_swe_ce}, we plot the dynamic social welfare equilibrium and the dynamic competitive equilibrium. From the plot, we can see that the dynamic social welfare equilibrium and the dynamic competitive equilibrium agree, which is consistent with Theorem \ref{thm4}. The dynamic optimal price for traded resource versus $t$ is also shown in Figure~\ref{fig:dym_price}, where the price  experiences oscillations both at the beginning and in the end over the time horizon, and holds a steady value in between. 

\section{Conclusions}\label{sec::conc}
We studied multi-agent systems with decentralized resource allocation without  external resource supply.  For multi-agent systems with static local allocations, we showed    that under general convexity assumptions, the competitive equilibrium and the social welfare equilibrium  exist and agree using a duality analysis. We also studied the problem of social shaping for competitive equilibriums, where   the pricing under a competitive equilibrium is associated with an upper bound. We presented an explicit family of socially admissible utility functions under which   the pricing   is always socially acceptable. Finally, a dynamical multi-agent system was studied and generalized in an optimal control context. In future work, numerical algorithms for the computation of  the socially admissible utility function are possible.


\begin{thebibliography}{1} 

\bibitem{tsi1984} J. N. Tsitsiklis, {\em Problems in decentralized decision making and computation}, Massachusetts Institute of Tech Cambridge Lab for Information \& Decision Systems, 1984. 

\bibitem{Jadbabaie03} A. Jadbabaie, J. Lin, and A. S. Morse, ``Coordination of groups of mobile autonomous agents using nearest neighbor rules," {\em IEEE Trans. Autom. Control}, vol. 48, no. 6, pp. 988–1001, 2003.

\bibitem{Murray04}  R. Olfati-Saber and R. M. Murray, ``Consensus problems in the networks of agents with switching topology and time delays," {\em IEEE Trans. Autom. Control},  vol. 49, no. 9, pp. 1520-1533, 2004.



\bibitem{Martinez07} S. Mart\'{i}nez, J. Cort\'{e}s, and F. Bullo, ``Motion coordination with distributed information," 
{\em IEEE Control Systems Magazine}, vol. 27, pp. 75–88, 2007. 



\bibitem{Nedic10} A. Nedi\'{c}, A. Ozdaglar, and P. A. Parrilo, “Constrained consensus and
optimization in multi-agent networks,” {\em IEEE Trans. Autom. Control,}
vol. 55, no. 4, pp. 922–938, 2010.

\bibitem{Mesbahi2010} M. Mesbahi and M. Egerstedt, {\em Graph Theoretic Methods in Multiagent Networks}, Princeton University Press, 2010.


\bibitem{Bert_N&P}
D. P. Bertsekas, \emph{Nonlinear Programming}, 2nd edition, Atehna Scientific,
  Belmont, Massachusetts, 1999.

\bibitem{boyd} S. Boyd and L. Vandenberghe, {\em Convex Optimization}, 1st edition, Cambridge University Press, 2004.

\bibitem{Acemoglu}  D. Acemoglu, D. Laibson, and J. List,  {\em Microeconomics}, 2nd edition, Pearson, 2018. 


\bibitem{Arrow1954} K. J. Arrow and G. Debreu, ``Existence of an equilibrium for a competitive economy," {\em Econometrica, Journal of the Econometric Society}, vol. 22, no. 3, pp. 265--290, 1954.

\bibitem{Debreu1952} G. Debreu, ``A social equilibrium existence theorem," {\em Proc. Nat. Acad. Sci.}, vol. 38, no. 10, pp. 886--893, 1952. 

\bibitem{CSM2020} S. Li, J. Lian, A. J. Conejo, and W. Zhang, ``Transactive energy systems: the market-based coordination of distributed energy resources,"  {\em IEEE Control Systems Magazine}, vol. 40, no. 4, pp. 26--52, 2020. 



\bibitem{Chen10} L. Chen, N. Li, S. H. Low, and J. C. Doyle, ``Two market models for demand response in power networks," in {\em Proc. IEEE SmartGridComm}, pp. 397--402, 2010.

\bibitem{Li2011} N. Li, L. Chen, and S. H. Low, ``Optimal demand response based on utility maximization in power networks," in {\em  Proc. IEEE Power and Energy Society General Meeting}, pp. 1–8, 2011. 

\bibitem{Papadaskalopoulos13}D. Papadaskalopoulos and G. Strbac, ``Decentralized participation of flexible demand in electricity marketsart i: Market mechanism," {\em IEEE Trans. Power Syst.}, vol. 28, no. 4, pp. 3658--3666, 2013.

\bibitem{Hansen15} J. Hansen, J. Knudsen, and A. Annaswamy, ``A dynamic market mechanism for integration of renewables and demand response," {\em IEEE Trans. Control Syst. Technol.}, vol. 24, no. 3, pp. 940--955, 2015. 

  

\bibitem{Li2015} N. Li, L. Chen, and M. A. Dahleh, ``Demand response using linear supply function bidding," {\em IEEE Trans. Smart Grid}, vol. 6, no. 4, pp. 1827--1838, 2015.



\bibitem{Jadhav2018} A. M. Jadhav, N. R. Patne, and J. M. Guerrero, ``A novel approach to neighborhood fair energy trading in a distribution network of multiple microgrid clusters," {\em IEEE Trans. Ind. Electron.}, vol. 66, no. 2, pp. 1520--1531, 2018. 


\bibitem{Muthirayan2019}  D. Muthirayan, D. Kalathil, K. Poolla, and P. Varaiya, ``Mechanism design for demand response programs," {\em IEEE Trans. Smart Grid}, vol. 11, no. 1, pp. 61--73, 2019.  

 
\bibitem{Alvarado2001} F. L. Alvarado, J. Meng, C. L. DeMarco, and W. S. Mota, ``Stability analysis of interconnected power systems coupled with market dynamics," {\em IEEE Trans. Power Syst.}, vol. 16, no. 4, pp. 695–701, 2001.


\bibitem{Singh2018} R. Singh, P. R. Kumar, and L. Xie, ``Decentralized control via dynamic stochastic prices: The independent system operator problem," {\em IEEE Trans. Autom. Control}, vol. 63, no. 10, pp. 3206--3220, 2018. 


\bibitem{Stegink2017} T. Stegink, C. D. Persis,   and A. V. D. Schaft, ``A unifying energy-based approach to stability of power grids with market dynamics," {\em IEEE Trans Autom. Control}, vol. 62, no. 6, pp. 2612--2622, 2017. 


\bibitem{Zhang2015} X. Zhang and A. Papachristodoulou, ``A real-time control framework for
smart power networks: Design methodology and stability," {\em Automatica},
vol. 58, pp. 43–50, 2015.

\bibitem{Jovic2010} A. Jokic, M. Lazar, and P. P. J. V. D. Bosch,``Price-based control of electrical power systems," in {\em Intelligent Infrastructures},   New York: Springer-Verlag, pp. 109–131, 2010. 




\bibitem{Jovic2009} A. Jokic, M. Lazar, and P. P. V. D. Bosch, ``On constrained steady-state regulation: Dynamic KKT controllers," {\em IEEE Trans. Autom. Control},
vol. 54, no. 9, pp. 2250–2254, 2009.



\bibitem{Nordhaus96}W. D. Nordhaus and Z. Yang, “A regional dynamic general-equilibrium model of alternative climate-change strategies,” {\em The American Economic Review}, vol. 86, no. 4, pp. 741–765, 1996.

\bibitem{Nordhaus2011} W. D. Nordhaus, “Estimates of the social cost of carbon: background and results from the rice-2011 model,” {\em National Bureau of Economic Research}, no. w17540, 2011.

\bibitem{Nordhaus2017}W. D. Nordhaus, “Revisiting the social cost of carbon,” {\em Proc. Nat. Acad. Sci.}, vol. 114, no. 7, pp. 1518–1523, 2017.

\bibitem{Chris2019} C. M. Kellett, S. R. Weller, T. Faulwasser, L. Grunec, and 
W. Semmler, ``Feedback, dynamics, and optimal control in climate economics," {\em Annual Reviews in Control}, vol. 47, pp, 7--20, 2019. 


\bibitem{Hallam} A.  Hallam, ``Competitive equilibrium and societal welfare," {\em Lecture Notes}, 2016. 

\bibitem{CNN} C. Morehouse, ``Texas PUC loses 2nd commissioner as Lt. Gov. presses ERCOT to correct \$16B overcharge," {\em www.utilitydive.com}, 2021. 

\bibitem{Kuyper} J. Kuyper, H. Schroeder, B.-O. Linn\'{e}r, ``The evolution of the UNFCCC," {\em Annual Review of Environment and Resources}, vol. 43, pp. 343--368, 2018.

 
 


\end{thebibliography}
\end{document}